%% file: JCP_2017_CF_v8.tex
\documentclass[amsmath,amssymb,aip,jcp,reprint]{revtex4-1}

\usepackage{graphicx}
\usepackage{dcolumn}
\usepackage{bm}
\usepackage{epsf}
\usepackage{epsfig}
\usepackage{amsmath}
\usepackage{latexsym}
\usepackage{amssymb}
\usepackage{subcaption}
\usepackage{widetable}
\usepackage{cleveref}
\usepackage{tikz,pgfplots}

\newcommand{\e}{\ensuremath{\mathrm{e}}}

\newcommand{\R}{\mathbb{R}}

\newcommand{\structurecomment}[1]{}

\newcommand{\T}{\tau}

\renewcommand{\P}{P}

\begin{document}
\title[Exponential propagators for the Schr\"odinger equation]{Exponential propagators for the Schr\"odinger equation with a time-dependent potential}
\author{Philipp Bader}
\email{bader@uji.es}
\affiliation{Departament de Matem\`atiques, Universitat Jaume I,
E-12071 Castell\'on, Spain.
}
\author{Sergio Blanes}
\email{serblaza@imm.upv.es.}
\affiliation{Instituto de Matem\'atica Multidisciplinar,
Universitat Polit\`ecnica de Val\`encia, E-46022 Valencia, Spain.}
\author{Nikita Kopylov}
\email{nikop1@upvnet.upv.es}
\affiliation{ 
Instituto de Matem\'atica Multidisciplinar,
Universitat Polit\`ecnica de Val\`encia, E-46022 Valencia, Spain.}
\date{\today}

\begin{abstract}
We consider the numerical integration of the Schr\"o\-din\-ger equation with a time-dependent Hamiltonian given as the sum of the kinetic energy and a time-dependent potential.
Commutator-free (CF) propagators are exponential propagators that have shown to be highly efficient for general time-dependent Hamiltonians.
We propose new CF propagators that are tailored for Hamiltonians of said structure, showing a considerably improved performance.
We obtain new fourth- and sixth-order CF propagators as well as a novel sixth-order propagator that incorporates a double commutator that only depends on coordinates, so this term can be considered as cost-free. 
The algorithms require the computation of the action of exponentials on a vector similarly to the well known exponential midpoint propagator, and this is carried out using the Lanczos method. 
We illustrate the performance of the new methods on several numerical examples.
\end{abstract}

\keywords{Schr\"odinger equation, time-dependent potential, propagators, commutator-free methods}
\maketitle

\section{Introduction} \label{sec.1}
We study the numerical integration of the time-dependent Schr\"odinger equation (SE) in units such that the Planck constant $\hbar = 1$:
\begin{equation} \label{Schr0}
i \frac{\partial}{\partial t} \psi (x,t) = \widehat{H}(t) \psi (x,t), \quad \psi(x,t_0)=\psi_0(x),
\end{equation}
where the Hamiltonian is given by
\begin{equation}\label{eq:hamiltonian_sum}
\widehat{H}(t) = \widehat{T} + \widehat{V}(x,t) = -\frac{\Delta}{2\mu} + \widehat{V}(x,t), 
\end{equation}
$x\in\R^d$, $\mu$ is the reduced mass, $\Delta$ is the Laplacian operator and $\widehat{V}(x,t)$ is an explicitly time-dependent potential.

Although postponed spatial discretization has been proposed\cite{bader14eaf,bader16emf}, the first step is usually discretizing eq. \eqref{Schr0} taking $ N $ grid points.
For this sake, various techniques can be used (see Ref. \onlinecite{lubich08fqt} and references therein).
Thus, one obtains a linear ordinary differential equation:
\begin{equation} \label{eq:schr_disc}
i \frac{d }{dt} u(t) = H(t)\, u(t), \quad u(t_0)=u_{0} \in \mathbb{C}^N,
\end{equation}
where $H(t)=T+V(t)$ is an $N\times N$ Hermitian matrix.

As in the existing literature, we assume the problem to be either periodic in space $x$ or to be considered as one due to solution $ u(t) $ and its time-derivatives vanishing at the boundaries.

Thereafter, the time integration interval $ [t_0, t_f] $ is divided in a number of sufficiently small timesteps of length $\tau$, 
and stepwise approximations $ u_k $ to the solution values $ u(t_k)$ are evaluated at the temporal grid points $t_k = t_0 + \tau k$ ($k=1,2,\ldots$).

Many algorithms for solving equation (\ref{eq:schr_disc}) exist indeed in the literature:
split-operator methods \cite{feit82sot},
Runge--Kutta \cite{tremblay04upa} and symplectic partitioned Runge--Kutta methods\cite{sanzserna96cni}, a combination of a 4th-order Magnus method with the Lanczos algorithm \cite{kormann08atp} and the so-called $(t,t')$ method \cite{peskin94tso}.
In the papers \cite{castro04pft,kormann08atp}, a profound analysis of these schemes was carried out.

There are two highly efficient families of propagators that have been adapted\cite{blanes17sta,blanes15aea,sanzserna96cni,gray94chs} to the explicitly time-dependent case \eqref{eq:schr_disc}.
Both make use of the fact that $V$ is diagonal in coordinate space and fast Fourier transform (FFT) algorithms $ \mathcal{F} $ can be used to diagonalize the kinetic energy operator $T=\mathcal{F}^{-1}T_{D}\mathcal{F}$.  


In the first family, the solution is approximated by the unitary split operator algorithms, i.e. by compositions of the form
\begin{equation} \label{td.2a}
\e^{-i \tau (T+V)} \approx 
\e^{-i b_m \T V_m} \e^{-i a_m \T T} 
\cdots c
\e^{-i b_1 \T V_1} \e^{-i a_1 \T T},
\end{equation} 
where $ V_j = V(t_k + c_j \tau)$ and $\{ a_i, b_i, c_i \}$ are appropriately chosen real coefficients \cite{blanes16aci,mclachlan02sm,neuhauser09otc,thalhammer08hoe,thalhammer12cao}. In this case, $m$ FFT calls (and their inverses, IFFT) per step are required. 
This class of methods shows high performance when $\|T\|$, $\|V\|$, $\|[T,V]\| = \|TV-VT\|$ are relatively small or their action's results have small norms, otherwise high errors can occur.

On the other hand, if $H(t)$ and its action $H(t)u(t)$ satisfy similar norm conditions as above, splitting $ u(t) = q(t)+ip(t)$ into real and imaginary parts can be used.
Eq. \eqref{eq:schr_disc} becomes
\begin{equation} \label{SchrDiscr2}
\frac{d }{dt} \left( \begin{array}{c}
q\\
p
\end{array} \right) = 
\left[
\left( \begin{array}{ccc}
0 & & H(t) \\
0 & & 0 \end{array} \right) 
+
\left( \begin{array}{ccc}
0 & & 0 \\
-H(t) & & 0 \end{array} \right) 
\right]
\left( \begin{array}{c}
q\\
p
\end{array} \right)
\end{equation}
that is similar to a set of $N$ coupled harmonic oscillators.
Therefore, splitting methods tailored for that type of problems can be employed \cite{blanes15aea,gray96sit,sanzserna96cni,kormann08atp,gray94chs}. 
The algorithm requires $2m$ real-to-complex direct and inverse FFTs at a similar cost as $m$ complex-to-complex FFTs, and is superior to the Chebyshev approximation for short and long time integrations\cite{blanes15aea}.
This class of methods is valid for any real Hamiltonian (not necessarily separable into kinetic and potential part). These are symplectic methods that are conjugate to unitary ones. 
They are easy to apply, have low computational cost, but are conditionally stable: if $ \| H \| $ grows, one should reduce the time step for stability reasons.

In the present work, we consider commutator-free (CF) methods that address problems in which the above-mentioned techniques are not appropriate.
Existing CF methods were, however, created for general Hamiltonians and do not take into account specific structure \eqref{eq:hamiltonian_sum}.

The goal of this work is to construct fourth- and sixth-order geometric numerical integrators, which preserve most qualitative properties of the solution \cite{blanes16aci,hairer06gni,sanz-serna94nhp}, by composition of maps $\P_m$ that approximate 
\begin{equation} \label{eq:krylov_approx}
u_{k+1}=\e^{-i \tau H} u_{k} \approx \P_m(-i \tau H)\, u_{k},
\end{equation}
where $ P_m $ is an approximation in the $ m $-dimensional Krylov--Lanczos subspace \cite{lubich08fqt,park86uqt,Saad1992}.

In essence, the three new methods utilize the 6\textsuperscript{th}-order Gauss--Legendre (GL) quadrature nodes (but other quadrature rules can be used \cite{blanes17sta,blanes17tao}):
\begin{equation}
\label{GL6cj}
c_1 = \frac12 - \frac{\sqrt{15}}{10}, \quad c_2 = \frac12, \quad
c_3 = \frac12 +\frac{\sqrt{15}}{10},
\end{equation}
to calculate linear combinations $ \bar{V}_i=\sum_{j=1}^3a_{i,j}V_j $, with $V_j=V(t_n+c_j\tau)$ with $a_{i,j}$ specific to each method:
\begin{itemize}
\item the 4\textsuperscript{th}-order method
\begin{equation}\label{CF4_2exp}
\Upsilon^{[4]}_{2} = \e^{-i \tau \bar V_4} \ \e^{-i \frac{\T}2 (T+\bar V_{3})} \ \e^{-i \frac{\T}2 (T+\bar V_{2})} \ \e^{-i \tau \bar V_1};
\end{equation}
\item the 6\textsuperscript{th}-order method that incorporates into \eqref{CF4_2exp} one term, $ \widetilde{V} $, containing derivatives of the potential:
\begin{equation} \label{approx.4m}
\Upsilon_{2}^{[6]} = \e^{-i \tau (\bar V_4+\tau^2 \widetilde V)} \ \e^{-i \frac{\T}2 (T+\bar V_{3})} \ \e^{-i \frac{\T}2 (T+\bar V_{2})} \ \e^{-i \tau (\bar V_1+\tau^2 \widetilde V)};
\end{equation}
\item the 6\textsuperscript{th}-order one, with no derivatives of $V$:
\begin{equation}\label{CF6_3exp}
\Upsilon_{3}^{[6]}= \e^{-i \tau \bar V_5}
\e^{-i \tau (a_2 T+\bar V_{4})} 
e^{-i \tau (a_3 T+\bar V_{3})}
\e^{-i \tau (a_2T+\bar V_{2})}
\e^{-i \tau \bar V_1}.
\end{equation}
\end{itemize}
Since the first and last exponentials in the schemes \eqref{CF4_2exp}-\eqref{CF6_3exp} are diagonal matrices whose computational cost can be neglected, these methods require, in practice, the calculation of 2, 2 and 3 operator exponentials, respectively, which is done by the Krylov polynomial approximation \eqref{eq:krylov_approx}.

The resulting schemes are 5/3 to 3 times more cost efficient than the existing CF methods thanks to simplified algebraic structure of the Hamiltonian written as the sum of a kinetic energy and a time-dependent potential and the low cost of exponentiating the latter.

In \cref{sec:methodology}, we describe the general methodology applied to obtaining commutator-free methods, existing specimens of this class and their limitations.
In \cref{sec:new_methods}, we perform manipulations that lead to explicit forms of new methods.
The article concludes with numerical examples.

\section{Methodology}\label{sec:methodology}
\subsection{Commutator-free Magnus integrators}
For a sufficiently small time step, the solution of \eqref{eq:schr_disc} can be expressed as the exponential of the Magnus expansion\cite{blanes09tme,magnus54ote}
 $ \Omega (t) $.
Computing the action of the $ p^{th} $-order truncation, $\exp{\Omega^{[p]}}$, on a vector can be resource-consuming due to the commutators involved in $\Omega^{[p]}$.

To circumvent these overheads, one can approximate exponentials by a product of simpler matrices that do not contain commutators, e.g., a $r$-exponential method of order $p$ takes the form
\begin{equation} \label{eq:commutator-free}
\mathrm{CF}^{[p]}_r 
= \prod_{j=1}^r \e^{-i \tau \widetilde H_j} 
= \e^{\Omega(\T)} + {\cal O}(\T^{p+1}), 
\end{equation}
where 
\[
\widetilde{ H_j} =\sum_{l=1}^ma_{j,l} H(t_k+c_l\T) 
= {a}_j T + \sum_{l=1}^ma_{j,l} V(t_k+c_l\T)
\]
with ${a}_j = \sum_{l=1}^ma_{j,l}$. 
There exist highly efficient methods in the literature up to order eight \cite{alvermann11hoc,blanes06fas,thalhammer06afo,auer18mio}.
However, those methods are not optimized for problems in which the Hamiltonian is of form \eqref{eq:hamiltonian_sum}.

\subsection{Algebra and optimisation}
In this work, we use the 6\textsuperscript{th}-order GL quadrature with nodes \eqref{GL6cj}.
Then, three linear combinations of $H(t_k + c_j \tau)=:T+V_j$
\begin{eqnarray} 
\alpha_1 &=& \T (T+V_2)={\cal O}(\tau), \nonumber \\ 
\alpha_2 &=& \tau \frac{\sqrt{15}}{3} (V_3 - V_1)={\cal O}(\tau^2), \label{eq:alphas} \\
\alpha_3 &=& \tau \frac{10}{3} (V_3 - 2 V_2 + V_1)={\cal O}(\tau^3), \nonumber 
\end{eqnarray}
suffice\cite{blanes17tao,blanes09tme,munthekaas99cia} to build methods up to order 6. 

Since $\alpha_2,\ \alpha_3$ are diagonal, and $[\alpha_2, \alpha_3]=0$, a 6\textsuperscript{th}-order approximation $\Omega^{[6]}$ of the Magnus expansion $ \Omega $ is given by
\begin{equation} \label{sixth}
\Omega^{[6]} =\alpha_{1} + \frac{1}{12} \alpha_{3} - \frac{1}{12} [12] +
\frac{1}{360} [113] - \frac{1}{240} [212] + \frac{1}{720} [1112],
\end{equation}
where $[ij\ldots kl]$ represents the nested commutator $\lbrack\alpha_{i},[\alpha_{j}, [\ldots ,[\alpha_{k},\alpha_{l}]\ldots]]]$.
Moreover, condition $[\alpha_2,\alpha_3]=0$ makes that the coefficients need to satisfy a reduced number of order conditions to reach high orders.

Moreover, the operator $[212]$ is also represented by a diagonal matrix whose elements are spatial derivatives of $ \widehat{V}(x,t) $.
This term can be combined with $ \alpha_{2} $ and $ \alpha_{3} $ to improve performance. 
Consequently, some $a_j$ can be zeroed and the computational cost (in terms of FFT calls) of corresponding exponentials in \eqref{eq:commutator-free} can be neglected, as shown in the following section.

\section{Integrators for~the~Schr\"{o}dinger equation}\label{sec:new_methods}
\subsection{Fourth-order methods}
An $m$-exponential method of order four that makes use of $\alpha_1, \alpha_2, \alpha_3$ given in \eqref{eq:alphas} must satisfy:
\begin{align*}
\Upsilon_m^{[4]} &=
\prod_{i=1}^m \exp \left(	x_{i,1} \alpha_{1} + x_{i,2} \alpha_{2} + x_{i,3} \alpha_{3}	\right) \\
&= \exp\left( \alpha_{1} + \frac{1}{12} \alpha_{3} - \frac{1}{12} [12] \right) + {\cal O}(\T^{5}).
\end{align*}
Time symmetry is preserved\cite{blanes06fas} in a CF method when 
\[
x_{m-i+1,j}=(-1)^{j+1}x_{i,j}, \quad i=1,2,\ldots
\]
We take $m=3$, $x_{1,1}=0$, and the result is the symmetric composition:
\begin{align*}
\Upsilon_3^{[4]} 
&=\exp \left( -x_{1,2} \alpha_{2} + x_{1,3} \alpha_{3} \right) \\
&\times \exp \left( x_{2,1} \alpha_{1} + x_{2,3} \alpha_{3} \right)\\
&\times \exp \left( x_{1,2} \alpha_{2} + x_{1,3} \alpha_{3} \right) \,
\end{align*}
where the first and the last exponents are diagonal, and their computational cost is similar to a scheme with only one exponential, like the midpoint exponential method.
This composition has 4 parameters to solve the order equations to assure 4\textsuperscript{th} order.
Therefore, one free parameter remains for optimization. 
The scheme that satisfies the condition for $[113]$ has the solution \cite{bader16sif}
\[
x_{1,2}=-\frac{1}{12},\quad x_{1,3}=\frac{1}{60},\quad
x_{2,1}=1,\quad x_{2,3}=\frac{1}{20} .
\]

New 4\textsuperscript{th}-order schemes with additional free parameters are given by composition
\begin{equation}\label{Magnus_42}
\begin{split}
\Upsilon_{2}^{[4]} 
&=\exp \left(-x_{1,2} \alpha_{2} + x_{1,3} \alpha_{3} \right) \\ 
&\times \exp \left(x_{2,1} \alpha_{1} - x_{2,2} \alpha_{2} + x_{2,3} \alpha_{3} \right) \\
&\times \exp \left(x_{2,1} \alpha_{1} + x_{2,2} \alpha_{2} + x_{2,3} \alpha_{3}
\right) \\
&\times \exp \left(x_{1,2} \alpha_{2} + x_{1,3} \alpha_{3} \right).
\end{split}
\end{equation}
With two free parameters to satisfy the order conditions associated with $[113]$ and $[1112]$, we get the unique solution 
\[
x_{1,2}=-x_{1,3}=-\frac{1}{60},\ x_{2,1}=\frac{1}{2},\
x_{2,2}=-\frac{2}{15},\ x_{2,3}=\frac{1}{40}.
\]
Using \eqref{eq:alphas}, we can transform coefficients $x_{i,j}$ to $a_{i,j}$ as follows:
\[
 a_{i,j} =\sum_{k=1}^3 x_{i,k}G_{k,j}
\]
with
\begin{equation}\label{matrixG}
 G = \left( 
	\begin{matrix}
	0 & 1 & 0\\
	-\frac{\sqrt{15}}{3} & 0 & \frac{\sqrt{15}}{3} \\
	\frac{10}{3} & -\frac{20}{3} & \frac{10}{3}
	\end{matrix}
	\right)\!,
\end{equation}
%
and the method in terms of linear combinations of $ T $ and $ V_j $ becomes
\begin{equation} \label{approx.4}
\Upsilon^{[4]}_{2} = \e^{-i \tau \bar V_4} \ \e^{-i \frac{\T}2 (T+\bar V_{3})} \ \e^{-i \frac{\T}2 (T+\bar V_{2})} \ \e^{-i \tau \bar V_1},
\end{equation}
where
\begin{equation} \label{V_CF42b}
\begin{array} {l}
\bar V_1= a_{1,1} V_1 + a_{1,2}V_2 + a_{1,3}V_3, \\
\bar V_2= a_{2,1} V_1 + a_{2,2}V_2 + a_{2,3}V_3, \\
\bar V_3= a_{2,3} V_1 + a_{2,2}V_2 + a_{2,1}V_3, \\
\bar V_4= a_{1,3} V_1 + a_{1,2}V_2 + a_{1,1}V_3, 
\end{array}
\end{equation}
and the coefficients are
\[
\begin{array} {lll}
a_{1,1}=\frac{10+\sqrt{15}}{180}, &
a_{1,2}=-\frac19, &
a_{1,3}=\frac{10-\sqrt{15}}{180}, \\
a_{2,1}=\frac{15+8\sqrt{15}}{90}, &
a_{2,2}=\frac23, &
a_{2,3}=\frac{15-8\sqrt{15}}{90}.
\end{array}
\]

All but one (i.e. $\left[212\right]$) order conditions at $\tau^5$ are satisfied in this symmetric method. The scheme can be written as
\[
\Upsilon_{2}^{[4]} =
\exp\left(\Omega^{[6]} - z [212]+{\cal O}(\T^7)\right),\quad z=\frac{1}{21600},
\]
that results in a highly optimized fourth-order method.

\subsection{Sixth-order methods}	
Let us examine the operators
\begin{equation*} 
\widehat \alpha_1 = \tau \left(\widehat{T}+\widehat{V}_2 (t)\right)\!,
\ \widehat \alpha_2 = \tau \frac{\sqrt{15}}{3} 
	\left(\widehat{V}_3(t) - \widehat{V}_1(t) \right).
\end{equation*}
The action of commutator $ [\widehat{\alpha}_2,[\widehat{ \alpha_1},\widehat{ \alpha}_2]] $ can be expressed through spatial derivatives $ \widehat{V}'(x,t_k+c_j\tau)\equiv \widehat{V}_j' $:
\[
[\widehat{212}]\psi(x,t) =
-\tau^3\frac{5}{3\mu} 
	\left(\widehat{V}_3' - \widehat{ V}_1' \right)^2	\psi(x,t),
\]
and the corresponding matrix representation $[212]$ is also diagonal in coordinate space.
If spatial derivatives are relatively simple to evaluate, $ [212] $ can be used in \eqref{Magnus_42} without significantly increasing its cost, similarly to the high order force gradient algorithms in Refs.~\onlinecite{koseleff94fcf,chin97sif,omelyan02}, yielding a 6\textsuperscript{th}-order method given by:
\begin{equation}\label{Magnus_62}
\begin{split}
\Upsilon_{2}^{[6]}
&=\exp \left(-x_{1,2} \alpha_{2} + x_{1,3} \alpha_{3} + 
y \left[\alpha_{2},\alpha_{1},\alpha_{2}\right] \right) \\ 
&\times \exp \left(x_{2,1} \alpha_{1} - x_{2,2} \alpha_{2} + x_{2,3} \alpha_{3}
\right) \\
&\times \exp \left(x_{2,1} \alpha_{1} + x_{2,2} \alpha_{2} + x_{2,3} \alpha_{3}
\right) \\
&\times \exp \left(x_{1,2} \alpha_{2} + x_{1,3} \alpha_{3} + y \left[\alpha_{2},\alpha_{1},\alpha_{2}\right]\right).
\end{split}
\end{equation}
It is solved by the same $x_{i,j}$ as in \eqref{Magnus_42}, and 
\[
y=\frac{z}{2}=\frac{1}{43200}.
\] 
The coefficients of this method coincide with the scheme\cite{bader16sif} built for the numerical integration of the Hill's equation.

With \eqref{eq:alphas} and \eqref{V_CF42b}, \eqref{Magnus_62} transforms to
\begin{equation} \label{approx.4m}
\Upsilon_{2}^{[6]} = \e^{-i \tau (\bar V_4+\tau^2 \widetilde V)} \ \e^{-i \frac{\T}2 (T+\bar V_{3})} \ \e^{-i \frac{\T}2 (T+\bar V_{2})} \ \e^{-i \tau (\bar V_1+\tau^2 \widetilde V)},
\end{equation}
and thanks to the diagonal form of $ V(t) $
\[
\widetilde V = -\frac{5 y}{3\mu}
\left(V'(t_k + c_3 \tau)- V'(t_k+c_1 \tau) \right)^2.
\]

Moreover, if $V(t)$ consists of a constant part $V_{(c)}$ and a time-dependent one, $V_{(t)}$, the modified potential $\widetilde{V}$ only requires spatial derivatives of $V_{(t)}$.
In particular, for an external field problem $V_{(t)}=f(t)V_{(f)}$ one has
\[
\widetilde{V}u_k =
-\frac{1}{\mu}\frac{1}{25920} 
(f_3-f_1)^2 V_{(f)}'^2u_k,
\]
with $f_j=f(t_k+c_j\T)$ that is easily computed for predefined external interaction functions.

Nonetheless, in some cases $ V'(t) $ or an appropriate approximation can be difficult to obtain,
hence we consider schemes that have more exponentials but retain a relatively low computational cost.

Although 6\textsuperscript{th}-order CF methods with four exponentials and no derivatives have been obtained\cite{alvermann11hoc},
those schemes showed poor performance.
Additional exponentials have to be incorporated into the scheme to improve its performance \cite{alvermann11hoc,blanes06fas}. 
We insert an exponential in the middle, and the composition contains 2 free parameters.
There is one free parameter among $x_{1,3},x_{2,3}$ and $x_{3,3}$ that multiplies $\alpha_3$. 
Since these variables only appear linearly in two of the order conditions, we take $x_{1,3}=0$. 
To reduce computational cost, we use one parameter to eliminate contributions from $\alpha_{1}$ in the outermost exponentials:	
\begin{align*}
\Upsilon_{3}^{[6]} &=
\begin{aligned}[t] 
&\exp \left(-x_{1,2} \alpha_{2} \right) 
 \times \exp \left(
x_{2,1} \alpha_{1} - x_{2,2} \alpha_{2} + x_{2,3} \alpha_{3}
\right) \\
 & \times \exp \left(
x_{3,1} \alpha_{1} + x_{3,3} \alpha_{3}
\right) \\
 &\times \exp \left(
x_{2,1} \alpha_{1} + x_{2,2} \alpha_{2} + x_{2,3} \alpha_{3}
\right) 
 \times \exp \left(
x_{1,2} \alpha_{2} \right) 
\end{aligned}\\
&= \e^{-i \tau \bar V_5}
\e^{-i \tau (a_2 T+\bar V_{4})} 
e^{-i \tau (a_3 T+\bar V_{3})}
\e^{-i \tau (a_2T+\bar V_{2})}
\e^{-i \tau \bar V_1},
\end{align*}
where
\begin{equation*}
\begin{array} {l}
\bar V_1= a_{1,1} V_1 + a_{1,2}V_2 + a_{1,3}V_3, \\
\bar V_2= a_{2,1} V_1 + a_{2,2}V_2 + a_{2,3}V_3, \\
\bar V_3= a_{3,1} V_1 + a_{3,2}V_2 + a_{3,3}V_3, \\
\bar V_4= a_{2,3} V_1 + a_{2,2}V_2 + a_{2,1}V_3, \\
\bar V_5= a_{1,3} V_1 + a_{1,2}V_2 + a_{1,1}V_3, 
\end{array}
\end{equation*}
and we take the solution that minimizes the sum of the coefficients $|x_{i,j}|$, 
\begin{equation*}
\begin{split}
a_{1, 1} &= 0.01994096265093610745,\\
a_{1, 2} &= 0,\\
a_{1, 3} &= -a_{1, 1},\\
a_{2, 1} &= 0.4882524910228221957,\\
a_{2, 2} &= -0.0046136830175630621,\\
a_{2, 3} &= 0.0834019108602182940;\\
a_{3, 1} &= -0.29387662410526271191,\\
a_{3, 2} &= 0.4536718104795705687,\\
a_{3, 3} &= a_{3, 1},\\
a_2 &= \sum_{i=1}^{3}a_{2,i} = 0.56704071886547742757,\\
a_3&=\sum_{i=1}^{3} a_{3,i} =1-2a_2\\
&=-0.13408143773095485515.
\end{split}
\end{equation*}

The analysis of this family of methods carried out in Ref.~\onlinecite{hofstatter17nso} shows that there exist no 6\textsuperscript{th}-order commutator-free methods with all coefficients $x_{i,1}=a_i$ greater than zero.

\section{Numerical illustrations}\label{s:NI}
In this section we illustrate the methods' behaviour and performance on academic one-dimensional examples.

First, let us describe the general setup used for numerical experiments in this paper.
We define the wave function $\psi(x,t)$ on a sufficiently large spatial domain $[x_0,x_N)$ to ensure that its value and its derivatives vanish.
This allows us to impose periodic boundary conditions $\psi(x_{0},t)=\psi(x_{N},t)$ and hence the use of FFTs.  We divide the interval into $N$ bins of length $\Delta x=(x_N-x_0)/N$, hence $x_k=x_0+k\Delta x$, $k=0,\ldots,N-1$. 
The discrete vector $u(t) \in \mathbb{C}^N$ from \eqref{eq:schr_disc} has components $u_k=(\Delta x)^{1/2}\psi(x_{k},t)$, and its norm $\|u(t)\|_2$ does not depend on $ t$. 

Under these assumptions, we consider the one-dimensional Schr\"odinger equation ($\hbar=1$):
\begin{equation}\label{Schr1}
\begin{split}
i \frac{\partial}{\partial t} \psi (x,t) &= \left(
-\frac{1}{2\mu} \frac{\partial^2}{\partial x^2} +
V(x) + f(t) x \right) \psi (x,t),\\
\psi(x,0)&=\psi_0(x).
\end{split}
\end{equation}
The time-dependent part $ f(t)=A \cos(\omega t) $ corresponds to an external laser field.

To check accuracy, we calculate the reference solution at the final time $ u(t_f) $ using a sufficiently small time step.
Then, we compute the solution with each method for various numbers of timesteps and measure the 2-norm of the errors at the final time.
The cost of a scheme is counted in units of complex-to-complex FFT--IFFT pairs required for calculating exponentials in each step.

In all the examples, exponentials are approximated by the Lanczos method according to the procedure described in the appendix.

Finally, we plot the errors for each timestep versus the corresponding cost in double logarithmic scale.

%
%
%

In our experiments, we compare with the commutator-free methods obtained in Ref.~\onlinecite{alvermann11hoc}: the 4\textsuperscript{th}-order method with 3 exponentials (CF4:3Opt in the source notation, $ \mathrm{CF}_{3Opt}^{[4]}$ here); and the 6\textsuperscript{th}-order scheme with five exponentials ($\mathrm{CF}\text{6:5} =:\mathrm{CF}_{5}^{[6]} $). The 6\textsuperscript{th}-order scheme with six exponentials ($\mathrm{CF}\text{6:6}$) showed worst performance in our numerical experiments and it has not been considered.
Although in Ref.~\onlinecite{alvermann11hoc}, the authors present other optimized 6\textsuperscript{th}-order methods that use the nodes of the 8\textsuperscript{th}-order GL quadrature rule, our schemes can also use higher order quadrature rules (as shown in the appendix) and for this reason we limit ourselves to the comparison of propagators that are similarly based on the three-point GL rule.
Thus, each exponential is written in terms of the linear combinations of $H_i=H(t_k+c_i\tau)$ with $c_i, \ i=1,2,3$ given in \eqref{GL6cj}:
\begin{widetext}\begin{equation}
	\begin{split}
	\mathrm{CF}_{5}^{[6]}&= \prod_{i=1}^5 \e^{-i\tau(a_{i,1}H_1+a_{i,2}H_2+a_{i,3}H_3)},\\
	a_{1,1}&= 0.203952578716323,\ 
	a_{1,2}= -0.059581898090478,\
	a_{1,3}= 0.015629319374155;\\
	a_{2,1}&= 0.133906069544898,\ 
	a_{2,2}= 0.314511533222506,\ 
	a_{2,3}= -0.060893550742092;\\
	a_{3,1}&= -0.014816639115506,\ 
	a_{3,2}= -0.065414825819611,\ 
	a_{3,3}= -0.014816639115506;\\
	a_{i,j}&=a_{6-i,4-j}, \ i=4,5, \ j=1,2,3.
	{\tiny }	\end{split}
	\label{eq:CF6Alv}
	\end{equation}
\end{widetext}

We also consider the Magnus-based exponential methods where a second-order method is given by the first term in the expansion
\begin{equation} \label{MidpointKrylov}
\Upsilon^{[2]} = \exp\left(-i \int_{t_k}^{t_k+\T} H(t) dt\right),	 
\end{equation} 
that provides a second order approximation to the solution.
Then, it is sufficient to approximate the integral by the midpoint quadrature rule of order two:
\begin{equation} \label{MidpointKrylov}
\Upsilon^{[2]}_1 = \e^{-i \T (T+V(t_k + \tau/2))}
\approx \P_m \left(-i \tau \left(T+V_{\frac{\tau}{2}}\right)\right) .
\end{equation} 
However, as it is the case for all the new methods presented in this work, the schemes can be used using higher-order quadrature rules. For example, we can approximate the integral using the 6\textsuperscript{th}-order GL rule \eqref{GL6cj} as for all methods tested in this work, so use the following averaged one-exponential method is considered:
\begin{equation} \label{MidpointImproved}
\Upsilon^{[2]}_{1,3} = \e^{-i \T (T+\frac1{18}(5V_1+8V_2+5V_3))}	 
\end{equation} 
with $V_i=V(t_k + c_i\tau)$. 
It is worth noting that, if the integral approximation contains the dominant error of the method, $\Upsilon^{[2]}_{1,3}$ should provide more accurate result without increasing the computational cost.

The salient point is that the new  $\Upsilon^{[4]}_2$ and $\Upsilon^{[6]}_2$ are minor modifications of two consecutive steps of \eqref{MidpointImproved} and, in most cases, their overall cost is similar because the higher order methods can be used with a time step twice larger. 


\paragraph*{The Walker--Preston model} 
is a simple model that represents adequately many typical applications and may serve as an indicative benchmark.

We take the Morse potential $V(x)=D \left(1-\e^{-\alpha x}\right)^2 $ for $ x \in [-0.8, 4.32) $, subdivided into $N=64$ and $N=128$.
The parameters chosen are: $\mu= 1745 \ a.u.$, $ D= 0.2251 \ a.u.$ and $\alpha= 1.1741 \ a.u.$$ $
The amplitude is $A= A_0=0.011025 \ a.u.$ and the frequency is $\omega= \omega_0=0.01787$. This corresponds to a standard example of diatomic HF molecule in a strong laser field \cite{walker77qvc1}.

The numerical experiments are repeated with reduced intensity and frequency $A= A_0/2$ and $\omega= \omega_0/2$, leaving time steps the same.

As initial condition, we take the ground state of the Morse potential:
$$ \phi (x) =
\sigma \exp \left( -\left(\gamma-\frac{1}{2}\right)\alpha x \right)
\exp \left( -\gamma \e^{-\alpha x} \right),
$$
with $ \gamma=2D/w_0$, $w_0=\alpha \sqrt{2D/\mu}$ and $\sigma$ is a normalizing constant.
We integrate on $[0,10\: t_p]$ with $t_p=2\pi/\omega$. 

In Figure~\ref{fig:walker64} shows the results for $N=64$ while Figure~\ref{fig:walker128} shows the results for $N=128$. We observe that both new 6\textsuperscript{th}-order methods show the best performance when high accuracies are desired.
The 4\textsuperscript{th}-order $ \Upsilon_{2}^{[4]} $ shows comparable performance when large time steps are used and do not require the evaluation of the spatial derivative of the potential. We observe that the new methods show better performance than the exponential midpoint for all accuracies of practical interest and have similar complexity for their implementation so, any user of the exponential midpoint method should consider one of the new methods for solving this class of problems.
\begin{figure}
	\begin{subfigure}[b]{0.45\textwidth}
		\caption{$\omega=\omega_0,\ A=A_0$}
		\input{EvsC_schr_64_001787.tikz}
	\end{subfigure}
	~ 
	\begin{subfigure}[b]{0.45\textwidth}
		\caption{$\omega=\omega_0/2,\ A=A_0/2$}
		\input{EvsC_schr_64_0008935.tikz}		
	\end{subfigure}
	\caption{Efficiency plots in logarithmic scale of the 4\textsuperscript{th}- and 6\textsuperscript{th}-order methods for the Walker--Preston model with $ N=64 $. The new methods are plotted with thicker lines.}\label{fig:walker64}
\end{figure}
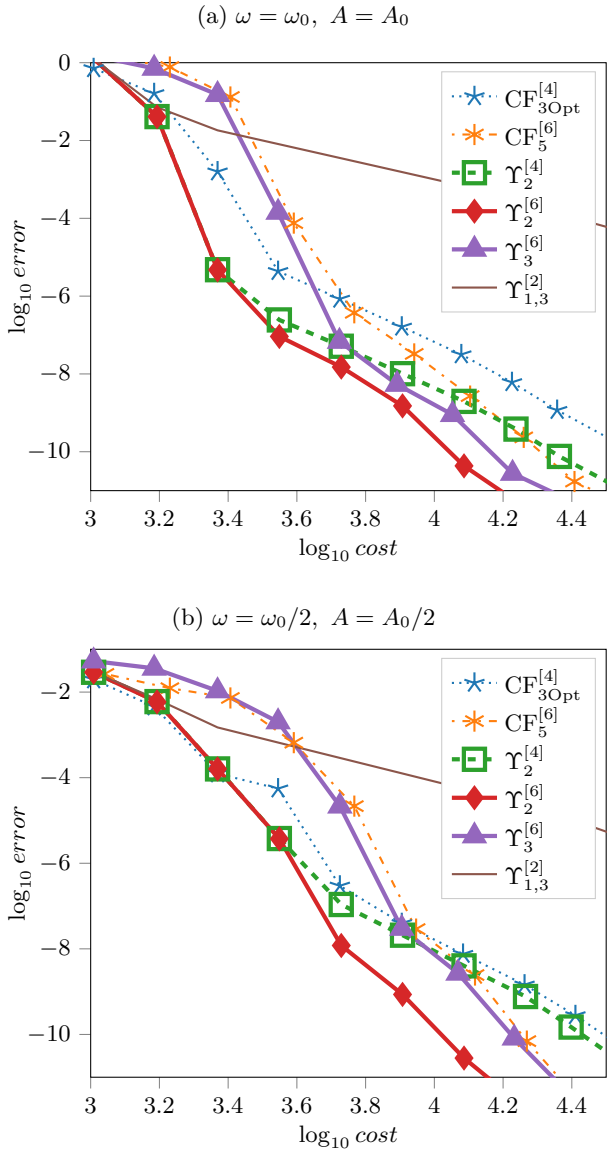
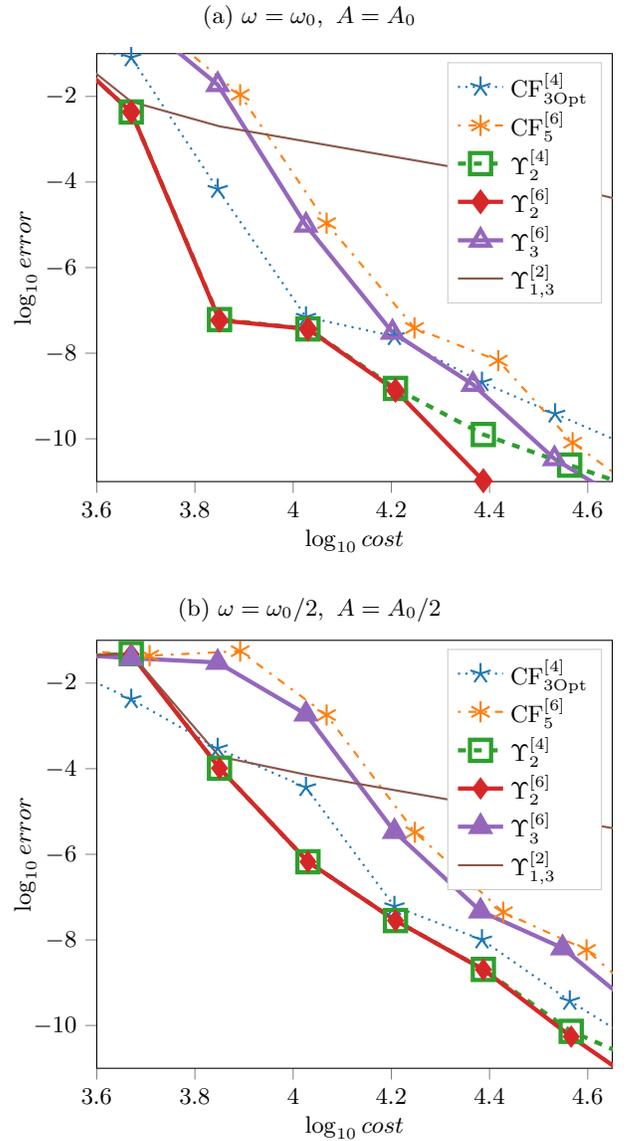
\begin{figure}
	\begin{subfigure}[b]{0.45\textwidth}
		\caption{$ \omega=\omega_0,\ A=A_0$}
		\input{EvsC_schr_128_001787.tikz}
	\end{subfigure}
	\begin{subfigure}[b]{0.45\textwidth}
		\caption{$  \omega=\omega_0/2,\ A=A_0/2$}
		\input{EvsC_schr_128_0008935.tikz}
	\end{subfigure}
	\caption{Same as Figure~\ref{fig:walker64} for $ N=128$.}\label{fig:walker128}
\end{figure}

\section{Conclusion}
We have constructed methods for the time-dependent Hamiltonian $ H(t)$ that can be presented as sum of the kinetic energy $ \widehat{T} $ and an explicitly time-dependent potential $ \widehat{V}(t) $.
The obtained methods belong to the class of so-called commutator-free (quasi-Magnus) integrators and heavily rely on the properties of the problem, hence their reduced computational cost.

On widely used benchmark problems, we have shown appropriateness of using the new methods to solve \eqref{Schr0}.
All new schemes have shown better performance than the non-customized methods of the same family.
Summarizing, $ \Upsilon_{2}^{[6]} $ is the method of choice when spatial derivatives of $ \widehat{V}(x,t) $ are available, and $ \Upsilon_{3}^{[6]} $ can be used otherwise.
The 4\textsuperscript{th}-order $ \Upsilon_{2}^{[4]} $ is appropriate for use with larger time steps when high accuracy is not required.
In addition, these methods have shown to be superior to the exponential midpoint for all desired accuracies of practical interest and have similar complexity for their implementation.

Future work may concern constructing of methods of higher order that employ other quadrature rules and can compete with 8\textsuperscript{th}-order CF methods\cite{alvermann11hoc} or the recent methods with simplified
commutators considered in Ref.~\onlinecite{iserles18mlm}.
Moreover, it seems possible to apply the same technique to tailoring methods for other problems with special structure which may allow additional optimization.

\appendix 
\section{Krylov approximation with Lanczos method}
If $ H $ is Hermitian and $ u $ is a unit vector, then $ \e^{-i \tau H}u_k $ can be approximated in the $ m $-dimensional Krylov subspace\cite{Saad1992}, whose column basis $ V_m=\left\lbrace v_1, v_2,\ldots,v_m \right\rbrace:=\left\lbrace u_k, Hu_k,\ldots,H^{m-1}u_k \right\rbrace $ is constructed by the Lanczos algorithm\cite{lubich08fqt}.
Exponentiation of the full operator is reduced to the exponentiation of a tridiagonal symmetric matrix $ T_m $ of smaller dimension:
\begin{equation}\label{eq:krylov_exp}
 \e^{-i \tau H}u_k \approx V_m \e^{-i \tau T_m} e_1, 
\end{equation}
where $ e_1 $ is the first column of the identity matrix.

The Lanczos algorithm builds an orthogonal basis and fills elements of $ T_m $\cite{lubich08fqt}:
\[
\begin{array}{l}
 {v}_{0} = {0}, \quad
 {v}_{1} = {u}_{k}, \quad
 \beta_{1} = {0}; \\
 {\bf for} \ \ i=1,\ldots, m \\
 \quad {y} = \tau H{v}_i - \beta_i {v}_{i-1} \\
 \qquad \alpha_i = \langle {v}_i | {y} \rangle \\
 \quad {y} = {y} - \alpha_i {v}_{i} \\
 \qquad \beta_{i+1} = \| {y} \| \\
 \quad {v}_{i+1} = {y}/\beta_{i+1} \\
 {\bf endfor} \\
 \alpha_i = (T_{m})_{i,i},\ 
  \beta_{i+1}=(T_{m})_{i+1,i},\ i=1,\ldots,m-1.
\end{array}
\]
The procedure requires $m$ matrix--vector products, storing $m$ resulting vectors, and at the end discards $v_{m+1}$ and $\beta_{m+1}$ to produce square matrices.

The error can be cheaply estimated by
\begin{multline*}
	err = \beta_{m+1}
	\left(\frac{2}{3} \left|e_m^T e^{-i \tau T_m/2} e_1\right|
	+\frac16\left|e_m^T e^{-i \tau T_m} e_1\right|\right),
	\end{multline*}
where $ e_1,\ e_m $ are the first and last column of the $m\times m$ identity matrix. 

Thus, in this work the procedure is the following: given a tolerance $ tol $, apply the Lanczos algorithm to build the sequence of vectors $\left\lbrace u_k, Hu_k, H^2 u_k,\ldots\right\rbrace $ and fill $ T_m $ until $ {err} < tol $ or $ m=10 $ is reached to ensure a sufficiently small error for a 6\textsuperscript{th}-order method.
Having obtained $ V_m$ and $ T_m $, then \eqref{eq:krylov_exp} is computed.

\section{Schemes with nodes of an arbitrary quadrature rule}
\label{otherqr}

The methods proposed in this work can be used with any other quadrature rule of order $r \ge 6$, say $\{\tilde b_i,\tilde c_i\}_{i=1}^k$, where $\tilde b_i$ are the weights and $\tilde c_i$ are the nodes. Similarly as in Ref.~\onlinecite{blanes17sta}, and as explicitly shown in Ref.~\onlinecite{blanes17tao}, there is a simple rule to apply: to replace $\alpha_1, \alpha_2,\alpha_3$ in \eqref{eq:alphas} by
\begin{equation*}
	\left( \begin{array}{c} 
\alpha_1 \\ \alpha_2 \\\alpha_3 
\end{array} \right)
= G \, Q^{-1} \tilde Q 
	\left( \begin{array}{c} 
\tilde H_1 \\ \vdots \\ \tilde H_k 
\end{array} \right)
\end{equation*}
with $\tilde H_i=H(t+\tilde c_i\tau), \ i=1,\ldots,k$, $G$ given by \eqref{matrixG} and
\[
 Q_{i,j}=b_j\left(c_j-\frac12\right)^{i-1}, \quad
 \tilde Q_{i,l}=\tilde b_l\left(\tilde c_k-\frac12\right)^{i-1}, 
\]
$i,j=1,2,3, \ l=1,\ldots, k$, with $(b_1,b_2,b_3)=\frac1{18}(5,8,5)$ being the weights of the 6\textsuperscript{th}-order GL quadrature rule with nodes $c_i$ given in \eqref{GL6cj}.

Obviously, the element $[212]$ to be used in the scheme $\Upsilon_{2}^{[6]}$ remains diagonal but now in terms the spatial derivatives of the potential evaluated at the new quadrature nodes.

\subsection*{Acknowledgments}
We wish to acknowledge Fernando Casas for its help in the construction of the methods $\Upsilon_{3}^{[6]}$.
The authors acknowledge Ministerio de Econom\'{\i}a y Competitividad (Spain) for financial support through project MTM2016-77660-P (AEI/FEDER, UE).
Additionally, Kopylov has been partly supported by grant GRISOLIA/2015/A/137 from the Generalitat Valenciana.

\bibliographystyle{aip}

\end{document}

%% file: EvsC_schr_64_001787.tikz
\begin{tikzpicture}

\definecolor{color0}{rgb}{0.12156862745098,0.466666666666667,0.705882352941177}
\definecolor{color1}{rgb}{1,0.498039215686275,0.0549019607843137}
\definecolor{color2}{rgb}{0.172549019607843,0.627450980392157,0.172549019607843}
\definecolor{color3}{rgb}{0.83921568627451,0.152941176470588,0.156862745098039}
\definecolor{color4}{rgb}{0.580392156862745,0.403921568627451,0.741176470588235}
\definecolor{color5}{rgb}{0.549019607843137,0.337254901960784,0.294117647058824}

\begin{axis}[
y label style={at={(axis description cs:0.08,0.5)},anchor=south},
xlabel={$\log_{10} cost$},
ylabel={$\log_{10} error$},
xmin=3, xmax=4.5,
ymin=-11, ymax=0,
tick align=outside,
tick pos=left,
x grid style={white!69.01960784313725!black},
y grid style={white!69.01960784313725!black},
legend style={draw=white!80.0!black},
legend cell align={left},
legend entries={{$\mathrm{ CF_{3Opt}^{[4]} }$},{$\mathrm{ CF_{5}^{[6]} }$},{$\Upsilon_{2}^{[4]}$},{$\Upsilon_{2}^{[6]}$},{$\Upsilon_{3}^{[6]}$},{$\Upsilon_{1,3}^{[2]}$}}
]
\addplot [thick, color0, dotted, mark=star, mark size=4, mark options={solid}]
table {%
2.47712125471966 0.150494007647374
2.65321251377534 0.176223311327777
2.81954393554187 0.188879894718167
3.00860017176192 -0.14539875734583
3.1846914308176 -0.787416348754147
3.36921585741014 -2.7991382853865
3.54530711646582 -5.35858026320271
3.72509452108147 -6.07128814103235
3.90525604874845 -6.78828987731045
4.07849306816657 -7.49999806418177
4.22587762863952 -8.21512527686531
4.35728673151412 -8.92991869418534
4.50490561894132 -9.64010201250091
4.65151060945951 -10.3375851192816
4.77558132418111 -10.9451235809566
};
\addplot [thick, color1, dash pattern=on 1pt off 3pt on 3pt off 3pt, mark=asterisk, mark size = 4, mark options={solid}]
table {%
2.69897000433602 0.0700283442406513
2.8750612633917 0.183032978498603
3.04139268515822 0.176200730496149
3.23044892137827 -0.111674705978656
3.40654018043396 -0.88108034014138
3.5910646070265 -4.12316870923992
3.7670073639498 -6.42877292974028
3.94146173934733 -7.48503542143241
4.10414555055401 -8.56908988589666
4.26011911003137 -9.63075790304711
4.40779866084079 -10.7658468757928
4.53992878789182 -11.3052228017522
4.68503375433325 -11.4901176426756
4.82230540982323 -11.3441934354434
4.98195444446997 -11.4334023931526
};
\addplot [ultra thick, color2, dashed, mark=square, mark size = 4, mark options={solid}]
table {%
2.30102999566398 0.115399698136327
2.47712125471966 0.192028113856127
2.64345267648619 0.166399815928079
2.83250891270624 0.176077437167733
3.00860017176192 0.102912189569723
3.19312459835446 -1.38549434955285
3.36921585741014 -5.32111415439712
3.54900326202579 -6.61175428374984
3.72916478969277 -7.28792756642655
3.90741136077459 -7.99402895196437
4.08635983067475 -8.70572766285299
4.2375688701982 -9.41773283449608
4.3638938977741 -10.1202475663122
4.51053160584908 -10.8067056268918
4.63978521298682 -11.2210731366329
};
\addplot [ultra thick, color3, mark=diamond*, mark size=3.5]
table {%
2.30102999566398 0.115399698136327
2.47712125471966 0.192084745157111
2.64345267648619 0.166382852402265
2.83250891270624 0.176072513187825
3.00860017176192 0.10291746030487
3.19312459835446 -1.38741224164016
3.36921585741014 -5.32527852009474
3.54900326202579 -7.03654857772259
3.72916478969277 -7.82363930247071
3.90741136077459 -8.82176351466955
4.08635983067475 -10.3634356962442
4.2375688701982 -11.2726394443022
4.3638938977741 -11.4619987598399
4.51053160584908 -11.3367133165421
4.63978521298682 -11.4070774961658
};
\addplot [ultra thick, color4, mark=triangle*, mark size=4]
table {%
2.47712125471966 0.174194073581594
2.65321251377534 0.076534845014341
2.81954393554187 0.16549889040036
3.00860017176192 0.178931124686568
3.1846914308176 -0.146267537692071
3.36921585741014 -0.820888349738911
3.54530711646582 -3.84795613560809
3.7238659644435 -7.16525322792732
3.89159320434897 -8.27050520038457
4.05292468370773 -9.04659059575938
4.2275782202996 -10.5580943177791
4.38334851777202 -11.2051086099927
4.51428204786038 -11.4494121252867
4.65420542968538 -11.3313926529635
4.79390208594707 -11.3887344217261
};
\addplot [thick, color5]
table {%
	2.30102999566398 0.116552653335919
	2.47712125471966 0.165874941348029
	2.64345267648619 0.169829536305536
	2.83250891270624 0.171633975953607
	3.00860017176192 0.0942858734522171
	3.19312459835446 -1.14486185086956
	3.36921585741014 -1.7385458502417
	3.54900326202579 -2.09728825066907
	3.72916478969277 -2.45725938074024
	3.90741136077459 -2.81360056516001
	4.08635983067475 -3.1714305609097
	4.23754373814287 -3.52926080708689
	4.36383752883614 -3.88613319884389
	4.51053160584908 -4.24350562036678
	4.63978521298682 -4.60064434548067
};

\end{axis}

\end{tikzpicture}

%% file: EvsC_schr_64_0008935.tikz
\begin{tikzpicture}

\definecolor{color0}{rgb}{0.12156862745098,0.466666666666667,0.705882352941177}
\definecolor{color1}{rgb}{1,0.498039215686275,0.0549019607843137}
\definecolor{color2}{rgb}{0.172549019607843,0.627450980392157,0.172549019607843}
\definecolor{color3}{rgb}{0.83921568627451,0.152941176470588,0.156862745098039}
\definecolor{color4}{rgb}{0.580392156862745,0.403921568627451,0.741176470588235}
\definecolor{color5}{rgb}{0.549019607843137,0.337254901960784,0.294117647058824}

\begin{axis}[
y label style={at={(axis description cs:0.08,0.5)},anchor=south},
xlabel={$\log_{10} cost$},
ylabel={$\log_{10} error$},
xmin=3., xmax=4.5,
ymin=-11, ymax=-1,
tick align=outside,
tick pos=left,
x grid style={white!69.01960784313725!black},
y grid style={white!69.01960784313725!black},
legend style={draw=white!80.0!black},
legend cell align={left},
legend entries={{$\mathrm{ CF_{3Opt}^{[4]} }$},{$\mathrm{ CF_{5}^{[6]} }$},{$\Upsilon_{2}^{[4]}$},{$\Upsilon_{2}^{[6]}$},{$\Upsilon_{3}^{[6]}$},{$\Upsilon_{1,3}^{[2]}$}}
]
\addplot [thick, color0, dotted, mark=star, mark size = 4, mark options={solid}]
table {%
2.47712125471966 -0.39348472755741
2.65321251377534 -0.228175529795905
2.81954393554187 -1.60373764371677
3.00860017176192 -1.70715751780248
3.1846914308176 -2.33939112534327
3.36921585741014 -3.90447601938435
3.54530711646582 -4.25602750821919
3.72509452108147 -6.52663201573848
3.90525604874845 -7.40249409856571
4.08350261983027 -8.1208408731148
4.26245108973043 -8.83816026689634
4.41049042000345 -9.55455111061459
4.54064225440965 -10.2686409473828
4.6638139591339 -10.983102930929
4.77871546654902 -11.3880146548704
};
\addplot [thick, color1, dash pattern=on 1pt off 3pt on 3pt off 3pt, mark=asterisk, mark size = 4, mark options={solid}]
table {%
2.69897000433602 0.0427708825367294
2.8750612633917 -0.743712300066745
3.04139268515822 -1.57519021003093
3.23044892137827 -1.90869595522815
3.40654018043396 -2.13704533408838
3.5910646070265 -3.18211165939467
3.76715586608218 -4.66462767187767
3.94694327069783 -7.53817595161424
4.1183308740573 -8.61167597902524
4.26913919679218 -10.1553382043235
4.41879829059035 -11.520886747325
4.57165061148795 -11.8940689032449
4.6990047465049 -11.7949757035965
4.85142958601948 -11.9124239655111
4.99132355943046 -12.8316693206281
};
\addplot [ultra thick, color2, dashed, mark=square, mark size = 4, mark options={solid}]
table {%
2.30102999566398 -0.100357857036747
2.47712125471966 -0.601785095650546
2.64345267648619 -1.18129597102487
2.83250891270624 -1.20331292915327
3.00860017176192 -1.53696005336369
3.19312459835446 -2.22260007483406
3.36921585741014 -3.7953988550214
3.54900326202579 -5.42249026789649
3.72916478969277 -6.96364242050047
3.90741136077459 -7.68070213487028
4.08635983067475 -8.39765563523231
4.26528962586083 -9.113770069104
4.39797475084092 -9.82775460455204
4.52551126096761 -10.5429626786744
4.64608990127443 -11.2290134127336
};
\addplot [ultra thick, color3, mark=diamond*, mark size=3.5]
table {%
2.30102999566398 -0.100357857036747
2.47712125471966 -0.640901273549984
2.64345267648619 -1.25084523178956
2.83250891270624 -1.20595940816361
3.00860017176192 -1.53699587672138
3.19312459835446 -2.22898361146428
3.36921585741014 -3.79548169307933
3.54900326202579 -5.42793621530183
3.72916478969277 -7.92160024337373
3.90741136077459 -9.06569112138169
4.08635983067475 -10.5523017576892
4.26528962586083 -11.7319269905014
4.39797475084092 -11.9224589089786
4.52551126096761 -11.9163334197485
4.64608990127443 -11.6709773692575
};
\addplot [ultra thick, color4, mark=triangle*, mark size=4]
table {%
2.47712125471966 0.0403928496809517
2.65321251377534 -0.20709645220323
2.81954393554187 -1.26800910994249
3.00860017176192 -1.27587521467098
3.1846914308176 -1.44459330856029
3.36921585741014 -1.9727663977175
3.54530711646582 -2.6973282584865
3.72509452108147 -4.6660142205449
3.90525604874845 -7.523437114606
4.06877936300956 -8.57156363122174
4.23248786635299 -10.0872098458617
4.39562339435584 -11.4207499255712
4.54288764234241 -11.8204847008823
4.6638139591339 -11.7357881699855
4.81698989530431 -11.7107742700248
};
\addplot [thick, color5]
table {%
	2.30102999566398 -0.899925957642784
	2.47712125471966 -0.155432733459199
	2.64345267648619 -1.16487172325727
	2.83250891270624 -1.37808471723275
	3.00860017176192 -1.52754807621379
	3.19312459835446 -2.16768429247866
	3.36921585741014 -2.8258345591548
	3.54900326202579 -3.18650065724928
	3.72916478969277 -3.54645785521204
	3.90741136077459 -3.90278874394003
	4.08635983067475 -4.26061396942508
	4.26528962586083 -4.618442079944
	4.39797475084092 -4.97531352187855
	4.52551126096761 -5.33268552220217
	4.64608990127443 -5.68982403480545
};

\end{axis}

\end{tikzpicture}

%% file: EvsC_schr_128_001787.tikz
\begin{tikzpicture}

\definecolor{color0}{rgb}{0.12156862745098,0.466666666666667,0.705882352941177}
\definecolor{color1}{rgb}{1,0.498039215686275,0.0549019607843137}
\definecolor{color2}{rgb}{0.172549019607843,0.627450980392157,0.172549019607843}
\definecolor{color3}{rgb}{0.83921568627451,0.152941176470588,0.156862745098039}
\definecolor{color4}{rgb}{0.580392156862745,0.403921568627451,0.741176470588235}
\definecolor{color5}{rgb}{0.549019607843137,0.337254901960784,0.294117647058824}

\begin{axis}[
y label style={at={(axis description cs:0.08,0.5)},anchor=south},
xlabel={$\log_{10} cost$},
ylabel={$\log_{10} error$},
xmin=3.6, xmax=4.65,
ymin=-11, ymax=-1,
tick align=outside,
tick pos=left,
x grid style={white!69.01960784313725!black},
y grid style={white!69.01960784313725!black},
legend style={draw=white!80.0!black},
legend cell align={left},
legend entries={{$\mathrm{ CF_{3Opt}^{[4]} }$},{$\mathrm{ CF_{5}^{[6]} }$},{$\Upsilon_{2}^{[4]}$},{$\Upsilon_{2}^{[6]}$},{$\Upsilon_{3}^{[6]}$},{$\Upsilon_{1,3}^{[2]}$}}
]
\addplot [thick, color0, dotted, mark=star, mark size = 4, mark options={solid}]
table {%
2.77815125038364 0.120509713598205
2.95424250943932 0.134553759956285
3.12057393120585 0.100258834986423
3.3096301674259 0.183052791023336
3.48572142648158 -0.673709129368067
3.67024585307412 -1.09158813740678
3.84633711212981 -4.1646628955265
4.02612451674545 -7.15374760696562
4.20628604441243 -7.60494754921118
4.38453261549425 -8.67891928895121
4.53351786201697 -9.41638760621319
4.67711398602634 -10.1300275664847
4.82468543313461 -10.7761917504629
4.95831537084576 -11.2485050339835
5.07970931648952 -11.3002416500918
};
\addplot [thick, color1, dash pattern=on 1pt off 3pt on 3pt off 3pt, mark=asterisk, mark size = 4, mark options={solid}]
table {%
3 0.0932203137425276
3.17609125905568 0.160462945052749
3.34242268082221 0.164886036732874
3.53147891704225 0.172565347588349
3.70757017609794 -0.210619461710691
3.89209460269048 -1.96213341924086
4.06818586174616 -4.96164250375397
4.24718739058635 -7.41049184273006
4.41783691162632 -8.17404567294395
4.56922170428736 -10.0911840219464
4.71982828625433 -11.3557099345873
4.8530651648216 -11.62118126504
4.9969536220331 -11.4788703125379
5.13996719518393 -11.4083980350386
5.28549879907422 -11.3811483661642
};
\addplot [ultra thick, color2, dashed, mark=square, mark size = 4, mark options={solid}]
table {%
2.60205999132796 0.149197141572398
2.77815125038364 0.107543591736859
2.94448267215017 0.120812006895421
3.13353890837022 0.133168405771965
3.3096301674259 0.177995320513486
3.49415459401844 -0.509502178148921
3.67024585307412 -2.36412991499708
3.85003325768977 -7.21441156589289
4.03019478535675 -7.43298075466265
4.20844135643857 -8.8234484686036
4.38738982633873 -9.89981319146018
4.56325557343525 -10.6155451167132
4.69201806227988 -11.1243238407964
4.82308932330194 -11.3772942099972
4.94711989693841 -11.3902363223389
};
\addplot [ultra thick, color3, mark=diamond*, mark size=3.5]
table {%
2.60205999132796 0.144362289691113
2.77815125038364 0.102080449696257
2.94448267215017 0.123912936887934
3.13353890837022 0.134803500640905
3.3096301674259 0.178074163250094
3.49415459401844 -0.509501981191282
3.67024585307412 -2.36412991253427
3.85003325768977 -7.23155359517211
4.03019478535675 -7.43471081131086
4.20844135643857 -8.86767943544469
4.38738982633873 -10.9805933232505
4.56325557343525 -11.6425362281474
4.69201806227988 -11.4694278358048
4.82308932330194 -11.4013668558788
4.94711989693841 -11.4033245762579
};
\addplot [ultra thick, color4, mark=triangle, mark size=4]
table {%
2.77815125038364 0.126225077895184
2.95424250943932 0.131456138141323
3.12057393120585 0.122432294479103
3.3096301674259 0.122275414943666
3.48572142648158 0.100718027953014
3.67024585307412 -0.0666245580296889
3.84633711212981 -1.72699451762121
4.02612451674545 -5.01073830661875
4.20237932107962 -7.51346302041444
4.36619874736973 -8.7320650047381
4.53223189257381 -10.468418540983
4.69665339001982 -11.6256681585865
4.8239434831132 -11.4451798995219
4.96470265802021 -11.3869908362344
5.11079502581682 -11.3842007982273
};
\addplot [thick, color5]
table {%
	2.60205999132796 0.08445394935806
	2.77815125038364 0.119007066968067
	2.94448267215017 0.127312892571154
	3.13353890837022 0.132439526288993
	3.3096301674259 0.177667714851941
	3.49415459401844 -0.506760723077473
	3.67024585307412 -2.12249666261456
	3.85003325768977 -2.69888427363071
	4.03019478535675 -3.05911628975629
	4.20844135643857 -3.41557115351346
	4.38738982633873 -3.7734513518138
	4.56325557343525 -4.13130361187914
	4.69201806227988 -4.48818569282119
	4.823095849973 -4.84556220967087
	4.94711989693841 -5.20270297975386
};
\end{axis}

\end{tikzpicture}

%% file: EvsC_schr_128_0008935.tikz
\begin{tikzpicture}

\definecolor{color0}{rgb}{0.12156862745098,0.466666666666667,0.705882352941177}
\definecolor{color1}{rgb}{1,0.498039215686275,0.0549019607843137}
\definecolor{color2}{rgb}{0.172549019607843,0.627450980392157,0.172549019607843}
\definecolor{color3}{rgb}{0.83921568627451,0.152941176470588,0.156862745098039}
\definecolor{color4}{rgb}{0.580392156862745,0.403921568627451,0.741176470588235}
\definecolor{color5}{rgb}{0.549019607843137,0.337254901960784,0.294117647058824}

\begin{axis}[
y label style={at={(axis description cs:0.08,0.5)},anchor=south},
xlabel={$\log_{10} cost$},
ylabel={$\log_{10} error$},
xmin=3.6, xmax=4.65,
ymin=-11, ymax=-1,
tick align=outside,
tick pos=left,
x grid style={white!69.01960784313725!black},
y grid style={white!69.01960784313725!black},
legend style={draw=white!80.0!black},
legend cell align={left},
legend entries={{$\mathrm{ CF_{3Opt}^{[4]} }$},{$\mathrm{ CF_{5}^{[6]} }$},{$\Upsilon_{2}^{[4]}$},{$\Upsilon_{2}^{[6]}$},{$\Upsilon_{3}^{[6]}$},{$\Upsilon_{1,3}^{[2]}$}}
]
\addplot [thick, color0, dotted, mark=star, mark size = 4, mark options={solid}]
table {%
2.77815125038364 -1.01155802479221
2.95424250943932 -0.946699640132459
3.12057393120585 -1.04824775332697
3.3096301674259 -1.30458875930572
3.48572142648158 -1.37230456108079
3.67024585307412 -2.3768099188373
3.84633711212981 -3.52403750842564
4.02612451674545 -4.43526697717548
4.20628604441243 -7.21785718360389
4.38453261549425 -7.99429154574715
4.56348108539441 -9.42661545643324
4.74241088058049 -10.703736069835
4.85870558937241 -11.3172174026807
5.02036128264771 -11.8960180737859
5.14341454208238 -11.8620628643015
};
\addplot [thick, color1, dash pattern=on 1pt off 3pt on 3pt off 3pt, mark=asterisk, mark size = 4, mark options={solid}]
table {%
3 -0.734918059693927
3.17609125905568 -1.04590977457271
3.34242268082221 -1.09426316945261
3.53147891704225 -1.20765501320529
3.70757017609794 -1.35781762986502
3.89209460269048 -1.25309856685217
4.06818586174616 -2.74560441681116
4.24797326636181 -5.48513981862287
4.42813479402879 -7.35396141260481
4.5976074408031 -8.23477412758426
4.74911766235632 -9.777168689491
4.89875808144042 -11.342812219734
5.02351583801841 -11.4895040669222
5.17872377474296 -11.5188027549582
5.29411776900455 -11.0784105320238
};
\addplot [ultra thick, color2, dashed, mark=square, mark size = 4, mark options={solid}]
table {%
2.60205999132796 -1.22108905904692
2.77815125038364 -0.907093030275879
2.94448267215017 -0.620525748395337
3.13353890837022 -1.1251387671684
3.3096301674259 -1.26713603867015
3.49415459401844 -1.4224419857577
3.67024585307412 -1.31060673410179
3.85003325768977 -3.98788671361352
4.03019478535675 -6.17420947505169
4.20844135643857 -7.54809853319258
4.38738982633873 -8.68963773495137
4.56631962152481 -10.137504217248
4.74476223706558 -11.0225665728337
4.87769377907898 -11.658503496407
5.0051118439161 -11.7699012866736
};
\addplot [ultra thick, color3, mark=diamond*, mark size =3]
table {%
2.60205999132796 -1.25950750369629
2.77815125038364 -0.894008404935139
2.94448267215017 -0.666784530385754
3.13353890837022 -1.12600212372472
3.3096301674259 -1.26555047924762
3.49415459401844 -1.42105378822211
3.67024585307412 -1.33134020492309
3.85003325768977 -3.98840436325719
4.03019478535675 -6.17423132599888
4.20844135643857 -7.54855651979485
4.38738982633873 -8.69290790042215
4.56631962152481 -10.2613691966974
4.74476223706558 -11.6800290154033
4.87769377907898 -11.872112638763
5.0051118439161 -11.7782859815244
};
\addplot [ultra thick, color4, mark=triangle*, mark size =4]
table {%
2.77815125038364 -1.07176219960852
2.95424250943932 -1.15844418919547
3.12057393120585 -1.07982672302727
3.3096301674259 -1.07578901715479
3.48572142648158 -1.28771726811236
3.67024585307412 -1.41757695788815
3.84633711212981 -1.5170236518829
4.02612451674545 -2.72528398260212
4.20628604441243 -5.47016557551895
4.38171058567103 -7.32881750726226
4.5487578285737 -8.18952548290252
4.71244765720305 -9.73020461072403
4.87509600556058 -11.3109368956546
5.00263251568728 -11.6161624232724
5.14341454208238 -11.4668064371432
};
\addplot [thick, color5]
table {%
	2.60205999132796 -0.663135757263645
	2.77815125038364 -0.865698457278351
	2.94448267215017 -0.916610203625236
	3.13353890837022 -1.08251141960164
	3.3096301674259 -1.24700428240419
	3.49415459401844 -1.42346373052151
	3.67024585307412 -1.31601617310297
	3.85003325768977 -3.7163114091566
	4.03019478535675 -4.1482803279515
	4.20844135643857 -4.50475260179833
	4.38738982633873 -4.86263188306041
	4.56631962152481 -5.22048361375193
	4.74476223706558 -5.57736542009512
	4.87769377907898 -5.93474190233759
	5.0051118439161 -6.29188247581978
};

\end{axis}

\end{tikzpicture}